\documentclass{amsart}
\usepackage[margin=1.5cm]{geometry}
\usepackage{hyperref}
\numberwithin{equation}{section}
\usepackage{amssymb}
\usepackage{amsmath, amsfonts,amsthm,amssymb,amscd, verbatim,graphicx,color,multirow,booktabs, caption,tikz,tikz-cd, mathdots,bm}
\usepackage{tikz-cd}
\usetikzlibrary{positioning}
\newtheorem{theorem}{Theorem}[section]
\newtheorem{example}{Example}[section]

\newtheorem{conjecture}[theorem]{Conjecture}

\def\nor#1#2{{\bf N}_{{#1}}({{#2}})}

\begin{document}

\title[Almost simple IBIS groups]{A classification of finite primitive IBIS groups with alternating socle}
\author{Melissa Lee}
\address{Department of Mathematics, University of Auckland, Auckland, New Zealand}
\email{melissa.lee@auckland.ac.nz}
\author{Pablo Spiga}
\address{Dipartimento di Matematica e Applicazioni, University of Milano-Bicocca, Via Cozzi 55, 20125 Milano, Italy} 
\email{pablo.spiga@unimib.it}

	\maketitle

        \begin{abstract}
   Let $G$ be a finite permutation group on $\Omega$. An ordered sequence $(\omega_1,\ldots,\omega_\ell)$ of elements of $\Omega$ is an irredundant base for $G$ if the pointwise stabilizer is trivial and no point is fixed by the stabilizer of its predecessors. If all irredundant bases of $G$ have the same cardinality, $G$ is said to be an IBIS group. 
   
Lucchini, Morigi and Moscatiello have proved a theorem reducing the problem of classifying finite primitive IBIS  groups $G$ to the case that the socle of $G$ is either abelian or non-abelian simple. 

In this paper, we classify the finite primitive IBIS groups having socle an alternating group. Moreover, we propose a conjecture aiming to give a classification of all almost simple primitive IBIS groups.
          \end{abstract}
\section{Introduction}
A \textit{\textbf{base}} for a permutation group $G$ is a sequence $(\omega_1,\ldots,\omega_\ell)$ whose pointwise stabilizer in $G$ is the identity. A base is said to be \textit{\textbf{irredundant}} if no base point is fixed by the stabilizer of the earlier points in the base. Therefore, this gives rise to a decreasing sequence of stabilizers
$$G>G_{\omega_1}>G_{\omega_1,\omega_2}>\cdots>G_{\omega_1,\omega_2,\ldots,\omega_{\ell-1}}>G_{\omega_1,\omega_2,\ldots,\omega_{\ell-1},\omega_\ell}=1.$$
The word redundant means that redundant points can be dropped and a base of smaller cardinality is obtained.

The minimal cardinality of a base is usually denoted by $b(G)$ and has played a key role in the investigation of primitive groups. (The number $b(G)$ is usually refereed to as the \textit{\textbf{base size}} of $G$.) For instance, if $G$ has a base of size $d$ and has degree $n$, then $|G|\le n^d$ and hence it is algorithmically very important to determine $b(G)$. There are many other statistics pertaining to primitive permutation groups and  related to the concept of base: relational complexity, height, maximal size of a minimal base, independent sets. All of these statistics are intertwined and have been investigated for instance in~\cite{GLS,MC}.

In this paper we are interested in the concept of IBIS group, this is an acronym for \textit{\textbf{Irredundant Bases of Invariant Size}}. Cameron and Fon-der-Flass~\cite{CF} (see also~\cite[Section~4.14]{Peter}) have proved that, in a finite permutation group the following conditions are equivalent:
\begin{itemize}
\item all irredundant bases have the same size;
\item the irredundant bases are invariant under re-ordering;
\item the irredundant bases are the bases of a matroid.
\end{itemize}
A permutation group satisfying one, and hence all, of these conditions is said to be IBIS. In particular, in an IBIS group $G$ all bases have cardinality $b(G)$ and this value is sometimes referred to as the \textit{\textbf{rank}} of the IBIS group. Except for this introductory section, we try not to use the term rank with this meaning, because it can be confusing with other concepts named rank in permutation groups.

There is a special circumstance where IBIS  groups arise. A transitive permutation group $G$ is said to be \textit{\textbf{geometric}} if $G$ permutes transitively its irredundant bases. Clearly, every geometric group is IBIS. Any regular permutation group is geometric with rank $1$. More remarkably,  Maund~\cite{Maund} has classified all finite geometric groups of rank greater than $1$.

Motivated by the classification of geometric groups, Cameron~\cite[page~124]{Peter} asks for a classification of the IBIS groups, specifying ``...it might be reasonable to pose this question for transitive IBIS groups (or perhaps primitive ones)".

Lucchini, Morigi and Moscatiello~\cite{LMM} have made the first attempt of classifying finite primitive IBIS groups. Their approach is via the O'Nan-Scott classification of primitive groups. 
\begin{theorem}[{Theorem~1.1,~\cite{LMM}}]Let $G$ be a primitive IBIS group. Then one of the following holds:
\begin{enumerate}
\item $G$ is of affine type,
\item $G$ is almost simple,
\item $G$ is of diagonal type.
\end{enumerate}
Moreover, $G$ is a primitive IBIS group of diagonal type if and only if it belongs to the infinite family of diagonal groups $\{\mathrm{PSL}_2(2^f)\times \mathrm{PSL}_2(2^f)\mid f\in\mathbb{N},f\ge 2\}$ having degree $|\mathrm{PSL}_2(2^f)|=2^f(4^f-1)$.
\end{theorem}
In the light of this result, the problem of understanding finite primitive IBIS groups is reduced to affine groups and to almost simple groups. Dealing with the first family seems very difficult. However, in our opinion, there is some hope for dealing with the second family. Indeed, suppose that $G$ is an almost simple primitive  IBIS group with $b(G)=2$. In particular, all bases of $G$ have cardinality $2$. Thus, for any two distinct points $\alpha$ and $\beta$ in the domain of $G$ we have $G_\alpha\cap G_\beta=1$. Therefore $G$ is a Frobenius group, contradicting the fact that $G$ is almost simple. This little argument shows that, if $G$ is an almost simple primitive IBIS group, then $b(G)\ge3$. Therefore, there is some hope in using the monumental work of Burness, O'Brien, Liebeck, Saxl, Shalev  and Wilson on the base size of primitive permutation groups, see for instance~\cite{burnessSym,burness1,burness2,burness3}.

In this paper, more modestly, we determine the almost simple primitive IBIS groups having socle an alternating group.
\begin{theorem}\label{main}Let $G$ be an almost simple primitive IBIS group having socle an alternating group. Then one of the following holds:
\begin{enumerate}
\item $G=\mathrm{Alt}(n)$ or $G=\mathrm{Sym}(n)$ of degree $n$;
\item $G=\mathrm{Alt}(5)=\mathrm{PSL}_2(5)$ or $G=\mathrm{Sym}(5)=\mathrm{PGL}_2(5)$ of degree $6$;
\item $G=\mathrm{Alt}(6)$ or $G=\mathrm{Sym}(6)$ of degree $6$ (on the right cosets of $\mathrm{PSL}_2(5)$ or $\mathrm{PGL}_2(5)$);
\item $G=\mathrm{Alt}(7)$  of degree $15$ (from the embedding $\mathrm{Alt}(7)<\mathrm{Alt}(8)=\mathrm{SL}_4(2)$; there are two such actions);
\item $G=\mathrm{Alt}(8)=\mathrm{SL}_4(2)$  of degree $15$ (there are two such actions);
\item $G=\mathrm{Alt}(6)$ or $G=\mathrm{Sym}(6)$ of degree $15$ acting on the $2$-subsets of $\{1,\ldots,6\}$;
\item $G\in \{\mathrm{Alt}(6)=\mathrm{PSL}_2(9),M_{10},\mathrm{P}\Sigma\mathrm{L}_2(9),\mathrm{P}\Gamma\mathrm{L}_2(9)\}$ of degree $10$. 
\end{enumerate}
\end{theorem}

Finally, in Section~\ref{examples}, we collect all the examples that I am aware of almost simple primitive IBIS groups. These examples naturally fall into ten families. As a wishful thinking we propose the following.

\begin{conjecture}\label{cong}
{\rm Let $G$ be an almost simple primitive IBIS group. Then $G$ is one of the groups appearing in Examples~\ref{example1}--\ref{example10} in Section~\ref{examples}.}
\end{conjecture}


\section{Almost simple groups with socle a non-abelian alternating group}\label{sec:alt}
In this section, we let $n$ be an integer $n\ge 5$ and we let $G$ be an almost simple group having socle the alternating group $\mathrm{Alt}(n)$. Recall that, except when $n=6$, we have
 $G=\mathrm{Alt}(n)$ or $G=\mathrm{Sym}(n)$. When $n=6$, the veracity of Theorem~\ref{main} follows from a computer computation. Therefore, for the rest of our argument, we let  $G\in\{\mathrm{Alt}(n),\mathrm{Sym(n)}\}$ act primitively and faithfully on a set $\Omega$ and we let $\omega\in \Omega$. We distinguish three cases, depending on whether $G_\omega$ is primitive, imprimitive or intransitive on $\{1,\ldots,n\}$.

\subsection{$G_\omega$ is primitive on $\{1,\ldots,n\}$.}
 From the main result in~\cite{burnessSym}, we see that either $b(G)=2$ or $n\le 12$. As $G$ is not a Frobenius group, $G$ cannot be IBIS when $b(G)=2$; hence, we deduce that $n\le 12$. Here the analysis is via a computer computation. When $n:=5$, we obtain the natural actions of $\mathrm{Alt}(5)=\mathrm{PSL}_2(5)$ and $\mathrm{Sym}(5)=\mathrm{PGL}_2(5)$ of degree $6$. When $n:=6$, we obtain the natural actions of $\mathrm{Alt}(6)$ and $\mathrm{Sym}(6)$ of degree $6$, on the right cosets of $\mathrm{PSL}_2(5)$ and $\mathrm{PGL}_2(5)$. When $n:=7$, we obtain the two  natural actions of $\mathrm{Alt}(7)$ of degree $15$ (we have two actions of this type and they arise from the embeddings of $\mathrm{Alt}(7)$ in $\mathrm{Alt}(8)=\mathrm{SL}_4(2)$, the stabilizer of a point is isomorphic to $\mathrm{SL}_3(2)$).  
When $n:=8$, we obtain the natural actions of $\mathrm{Alt}(8)=\mathrm{SL}_4(2)$ of degree $15$ (there are two such actions, on points and on hyperplanes of the $3$-dimensional projective space over the field with $2$ elements).

When $9\le n\le 12$, no more examples arise.

\subsection{$G_\omega$ is intransitive on $\{1,\ldots,n\}$.} The maximality of $G_\omega$ in $G$ gives that $G_\omega$ has two orbits on $\{1,\ldots,n\}$. Let $k$ and $n-k$ be the cardinality of these two orbits and, without loss of generality, assume $k<n/2$. Thus, we may identity the action of $G$ on $\Omega$ with the natural action of $G$ on the subsets of $\{1,\ldots,n\}$ having cardinality $k$. When $k=1$, we recover the natural action of $G$ of degree $n$, which is IBIS. Therefore, suppose $k>1$.

For the time being, suppose $k\ge 3$. In particular, $n>7$. Now, consider the following collection of $k$-subsets of $\{1,\ldots,n\}$:
$$\alpha_1:=\{1,\ldots,k-1,k\},\,\alpha_2:=\{1,\ldots,k-1,k+1\},\,\alpha_3:=\{1,\ldots,k-1,k+2\},\ldots, \alpha_{n-k}:=\{1,\ldots,k-1,n-1\}.$$
Clearly, we have a stabilizer chain
$$G_{\alpha_1}>G_{\alpha_1,\alpha_2}>G_{\alpha_1,\alpha_2,\alpha_3}>\cdots >G_{\alpha_1,\alpha_2,\ldots,\alpha_{n-k}}=G\cap \mathrm{Sym}(\{1,\ldots,k-1\}).$$
Set $H:=G\cap \mathrm{Sym}(\{1,\ldots,k-1\})$. Now, when $G=\mathrm{Alt}(n)$, consider the following collection of $k$-subsets of $\{1,\ldots,n\}$:
$$\beta_1:=\{1,k,k+1,\ldots,2k-2\},\,\beta_2:=\{2,k,k+1,\ldots,2k-2\},\,\ldots, \beta_{k-3}:=\{k-3,k,k+1,\ldots,2k-2\}.$$
Similarly,
when $G=\mathrm{Sym}(n)$, consider also the following $k$-subset:
$$\beta_{k-2}:=\{k-2,k,k+1,\ldots,2k-2\}.$$
Clearly, we have a stabilizer chain
$$H>H_{\beta_1}>H_{\beta_1,\beta_2}>\cdots >H_{\beta_1,\beta_2,\ldots,\beta_{k-3}}=H\cap \mathrm{Sym}(\{k-2,k-1\}).$$
When $G=\mathrm{Alt}(n)$, we obtain that $G$ has a base of cardinality $n-k+(k-3)=n-3$. Similarly, when $G=\mathrm{Sym}(n)$, by using also the $k$-subset $\beta_{k-2}$, we see that $G$ has a base of cardinality $n-k+(k-2)=n-2$.

Let us denote by $\kappa_3$ and $\kappa_k$ the base size of $\mathrm{Sym}(n)$ in its action on $2$-subsets and on $k$-subsets, respectively. From the work in~\cite{caceres,halasi}, we see that $$\kappa_k\le \kappa_3=\left\lceil\frac{n-1}{2}\right\rceil,$$
for $n\ge 9$.

Suppose now $G=\mathrm{Sym}(n)$.  From what we have seen above, in the action of $G$ on $k$-subsets of $\{1,\ldots,n\}$, there are bases having cardinality $n-2$ and bases having cardinality $\kappa_k\le \kappa_3$. Therefore, for $n \ge 9$, if $G$ is IBIS, we must have
$$n-2\le\left\lceil\frac{n-1}{2}\right\rceil. 
$$
However, this inequality never holds. It remains to consider the case $n< 9$. As $n>7$, we have $n=8$ and $k=3$. In this case, it can be verified that there are bases of cardinality $4$ and $6$: for instance
\begin{align*}
&(\{1,2,3\},\,\{1,2,4\},\,\{1,2,5\},\,\{1,2,6\},\,\{1,2,7\},\,\{1,3,4\}),\\
&(\{1,2,3\},\,\{4,5,6\},\,\{1,4,7\},\,\{2,5,8\}).
\end{align*}
In particular, $\mathrm{Sym}(8)$ in its action on $3$-subsets is not IBIS.

Suppose now $G=\mathrm{Alt}(n)$.  From what we have seen above, in the action of $G$ on $k$-subsets of $\{1,\ldots,n\}$ there are bases having cardinality $n-3$ and bases having cardinality at  most $\kappa_k\le \kappa_3$. Therefore, for $n \ge 9$, if $G$ is IBIS, we must have
$$n-3\le\left\lceil\frac{n-1}{2}\right\rceil. 
$$
However, this inequality never holds. It remains to consider the case $n< 9$. As $n>7$, we have $n=8$ and $k=3$. In this case, it can be verified that there are bases of cardinality $3$ and $4$: for instance
\begin{align*}
&(\{1,2,3\},\,\{4,5,6\},\,\{1,4,7\},\,\{2,5,8\}),\\
&(\{1,2,3\},\,\{1,4,5\},\,\{2,3,6\}).
\end{align*}
In particular, $\mathrm{Alt}(8)$ in its action on $3$-subsets is not IBIS.

We remains to deal with the case that $k=2$. In particular, $n\ge 5$. Now, consider the following collection of $2$-subsets of $\{1,\ldots,n\}$:
$$\alpha_1:=\{1,2\},\,\alpha_2:=\{1,3\},\,\alpha_3:=\{1,4\},\ldots \alpha_{n-3}:=\{1,n-2\}.$$
Clearly, we have a stabilizer chain
$$G_{\alpha_1}>G_{\alpha_1,\alpha_2}>G_{\alpha_1,\alpha_2,\alpha_3}>\cdots >G_{\alpha_1,\alpha_2,\ldots,\alpha_{n-3}}=G\cap \mathrm{Sym}(\{n-1,n\}).$$
When $G=\mathrm{Alt}(n)$, we obtain that $G$ has a base of cardinality $n-3$. Similarly, when $G=\mathrm{Sym}(n)$, by using also the point $\{1,n\}$, we see that $G$ has a base of cardinality $n-2$.

Let us denote by $\kappa_2$ the base size of $\mathrm{Sym}(n)$ in its action on $2$-subsets. From the work in~\cite{caceres,halasi}, we see that $$\kappa_2=\left\lceil\frac{2(n-1)}{3}\right\rceil.$$

Suppose now $G=\mathrm{Sym}(n)$.  From what we have seen above, in the action of $G$ on $2$-subsets of $\{1,\ldots,n\}$, there are bases having cardinality $n-2$ and bases having cardinality $\kappa_2$. Therefore,  if $G$ is IBIS, we must have
$$n-2=\left\lceil\frac{2(n-1)}{3}\right\rceil. 
$$
However, this equality  holds only when $n=6$. 

Suppose now $G=\mathrm{Alt}(n)$.  From what we have seen above, in the action of $G$ on $2$-subsets of $\{1,\ldots,n\}$, there are bases having cardinality $n-3$ and bases having cardinality at  most $\kappa_2$. Therefore, if $G$ is IBIS, we must have
$$n-3\le\left\lceil\frac{2(n-1)}{3}\right\rceil. 
$$
However, this inequality  holds only when $n\in \{5,6,7,8,9\}$. 
With a computer computation, we see that the only example arising is $\mathrm{Alt}(6)$ of degree $15$.






\subsection{$G_\omega$ is imprimitive on $\{1,\ldots,n\}$.} The maximality of $G_\omega$ in $G$ gives that $G_\omega$ leaves invariant only one system of imprimitivity $\Sigma$ for the action on $\{1,\ldots,n\}$. Let $b$ be the number of parts in $\Sigma$ and  let $a$ be the cardinality of the parts in $\Sigma$. Thus $n=ab$. Thus, we may idenfity the action of $G$ on $\Omega$ with the natural action of $G$ on the partitions of $\{1,\ldots,n\}$ having $b$ parts each having cardinality $a$. 
Suppose first that $a=2$. When $b\ge 4$, we have that the base size of $\mathrm{Sym}(2b)$ and $\mathrm{Alt}(2b)$ is $3$ in both cases, see for instance~\cite{MS}. Now consider the following $(2,b)$-regular partitions:
\begin{align*}
\alpha_1:=\{\{1,2\},\{3,4\},\{5,6\},\{7,8\},\ldots,\{2b-1,2b\}\},\\
\alpha_2:=\{\{1,4\},\{3,2\},\{5,6\},\{7,8\},\ldots,\{2b-1,2b\}\},\\
\alpha_3:=\{\{1,6\},\{3,4\},\{5,2\},\{7,8\},\ldots,\{2b-1,2b\}\},\\
\alpha_4:=\{\{1,8\},\{3,4\},\{5,6\},\{7,8\},\ldots,\{2b-1,2b\}\}.
\end{align*}
It is easy to verify that $G>G_{\alpha_1}>G_{\alpha_1,\alpha_2}>G_{\alpha_1,\alpha_2,\alpha_3}>G_{\alpha_1,\alpha_2,\alpha_3,\alpha_4}$
and hence $G$ is not IBIS. When $b=3$, it can be verified that $\mathrm{Alt}(6)$ and $\mathrm{Sym}(6)$ are both IBIS of degree $15$.

Assume now $a>2$.
Now consider the following $(2,b)$-regular partitions:
\begin{align*}
\alpha_1&:=\{\{1,2,3,\ldots,a-1,a\},\{a+1,a+2,a+3,\ldots,2a-1,2a\},\{3a+1,3a+2,3a+3,\ldots,4a-1,4a\},\ldots,\\
&\qquad\{a(b-1)+1,\ldots,ab-1,ab\}\},\\
\alpha_2&:=\{\{a+1,2,3,\ldots,a-1,a\},\{1,a+2,a+3,\ldots,2a-1,2a\},\{3a+1,3a+2,3a+3,\ldots,4a-1,4a\},\ldots,\\
&\qquad\{a(b-1)+1,\ldots,ab-1,ab\}\},\\
\alpha_3&:=\{\{1,a+2,3,\ldots,a-1,a\},\{a+1,2,a+3,\ldots,2a-1,2a\},\{3a+1,3a+2,3a+3,\ldots,4a-1,4a\},\ldots,\\
&\qquad\{a(b-1)+1,\ldots,ab-1,ab\}\},\\
\cdots&\cdots\\
\alpha_a&:=\{\{1,2,3,\ldots,2a-1,a\},\{a+1,a+2,a+3,\ldots,a-1,2a\},\{3a+1,3a+2,3a+3,\ldots,4a-1,4a\},\ldots,\\
&\qquad\{a(b-1)+1,\ldots,ab-1,ab\}\}.
\end{align*}
It is not hard to verify that $$G>G_{\alpha_1}>G_{\alpha_1,\alpha_2}>\cdots>\bigcap_{i=1}^aG_{\alpha_i}$$
and working in a similar fashion we may build a stabilizer chain of length
$\lfloor b/2\rfloor a$. From~\cite{MS}, we see that, except for an handful of cases (namely $(a,b)\in \{(4,2),(3,6),(3,7),(4,7),(7,3)\}$), the base size of $G$ on $(a,b)$-regular partitions is at most
$\lceil\log_b(a+3)\rceil+1$. These exceptional cases can be checked with the help of a computer. In the remaining cases, if $G$ is IBIS, then
$$\lceil\log_b(a+3)\rceil+1\ge \left\lfloor \frac{b}{2}\right\rfloor a.$$
However, this happens only when $(a,b)\in \{(3,2),(3,3),(4,2)\}$. By checking these cases with a computer we see that only $\mathrm{Alt}(6)$ and $\mathrm{Sym}(6)$ of degree $10$ arise here.


\section{Known examples of almost simple IBIS groups}\label{examples}
In this section, we collect all known examples of finite almost simple primitive IBIS groups. We believe that there are no more.
\begin{example}\label{example1}
{\rm Let $G$ be an almost simple group with socle $\mathrm{PSL}_2(q)$ acting on the $q+1$ points of the projective line, where $q:=p^f$, for some prime number $p$ and for some positive integer $f$. In particular, $\mathrm{PSL}_2(q)\unlhd G\le\mathrm{P}\Gamma\mathrm{L}_2(q)$. Set
$$r:=|G\mathrm{PGL}_2(q):\mathrm{PGL}_2(q)|.$$ Then $G$ is IBIS if and only if
$$r=1, \hbox { or }r \textrm{ is prime and }r\mid f.$$

We give a sketch of the proof of this ``if and only if''. We denote the elements of the projective line with square brackets $[a,b]$, with $a,b\in\mathbb{F}_q$. The pointwise stabilizer, $H$ say, of the frame $[1,0],[0,1],[1,1]$ consists of the field automorphisms of $\mathbb{F}_{p^f}/\mathbb{F}_p$. Moreover, the action of $H$ on the remaining points is permutation equivalent to the Galois action of $H$ on $\mathbb{F}_{p^f}$. Let $K$ be the subfield of $\mathbb{F}_{p^f}$ fixed by $H$. Then $K=\mathbb{F}_{p^\ell}$ for some divisor $\ell$ of $f$ and $H$ is the Galois group of $\mathbb{F}_{p^f}/\mathbb{F}_{p^\ell}$. This action is IBIS if and only if the Galois group of $\mathbb{F}_{p^f}/\mathbb{F}_{p^\ell}$ has all orbits of the same size on $\mathbb{F}_{p^f}\setminus\mathbb{F}_{p^\ell}$. This happens only when $f=\ell$ (that is, $H=1$) or when $f/\ell=|H|$ is a prime number.
}
\end{example}
\begin{example}\label{example2}
{\rm The Mathieu groups $M_{11}, M_{12}, M_{22}, M_{23}$ and $M_{24}$ in their natural action are all IBIS.}
\end{example}

\begin{example}\label{example3}
{\rm Let $G$ be the group $\mathrm{PGL}_n(2)=\mathrm{SL}_n(2)$ in its primitive action of degree $2^n-1$, on the non-zero vectors of the $n$-dimensional vector space $\mathbb{F}_2^n$ over the finite field $\mathbb{F}_2$ having cardinality $2$.

Since a base for $G$ is simply a basis of the vector space $\mathbb{F}_2^n$, we deduce that $G$ is IBIS. }
\end{example}

\begin{example}\label{example4}
{\rm Let $G$ be the group $\mathrm{PSp}_n(2)=\mathrm{Sp}_n(2)$ in its primitive action of degree $2^n-1$, on the non-zero vectors of the $n$-dimensional vector space $\mathbb{F}_2^n$ over the finite field $\mathbb{F}_2$ having cardinality $2$.

Since a base for $G$ is simply a basis of the vector space $\mathbb{F}_2^n$, we deduce that $G$ is IBIS. 

When $n=4$, we have $\mathrm{Sp}_4(2)\cong\mathrm{Sym}(6)$ and it can be checked that also $\mathrm{Sp}_4(2)'\cong\mathrm{Alt}(6)$ is IBIS. }
\end{example}

\begin{example}\label{example5}
{\rm Let $G$ be the Suzuki group $\mathrm{Sz}(q)={}^2B_2(q)$ in its $2$-transitive action of degree $q^2+1$ on the points of the Suzuki ovoid $\Omega$. Here we may identify $\Omega$ with the elements of the set
$$\{(\eta_1,\eta_2,\eta_3)\in\mathbb{F}_q^3\mid \eta_3=\eta_1\eta_2+\eta_1^{\sigma+2}+\eta_2^\sigma\}\cup\{\infty\},$$
where $\sigma:\mathbb{F}_q\to \mathbb{F}_q$ is the field automorphism such that $\sigma^2:\xi\mapsto \xi^2$ $\forall\xi\in \mathbb{F}_q$.

Following~\cite{DM}, the stabilizer $H$ in $G$ of the points $\infty$ and $(0,0,0)$ consists of the permutations $\eta_\kappa$, with $\kappa\in \mathbb{F}_q\setminus\{0\}$, satisfying
$$(\eta_1,\eta_2,\eta_3)^{\eta_\kappa}=(\kappa \eta_1,\kappa^{\sigma+1}\eta_2,\kappa^{\sigma+2}\eta_3),$$
for each $(\eta_1,\eta_2,\eta_3)\in \Omega$. Therefore, the orbits of $H$ on $\Omega\setminus\{\infty,(0,0,0)\}$ all have the same cardinality and this cardinality equals $|H|=q-1$. In particular, $G$ is IBIS.

Observe that, if $G<X\le\mathrm{Aut}(G)$, then $X$ is not IBIS in its action on $\Omega$. Indeed, $(1,1,1)\in\Omega$ and the pointwise stabilizer in $X$ of $\infty,(0,0,0),(1,1,1)$ has cardinality $|X:G|$ and consists of the field automorphisms of $\mathbb{F}_q$ contained in $X$. Therefore, $X$ admits a base of cardinality at least $4$. Whereas, if we choose $(1,\eta_2,1+\eta_2+\eta_2^\sigma)\in \Omega$ with $\eta_2$ a field generator of $\mathbb{F}_q$, then the pointwise stabilizer in $X$ of $\infty,(0,0,0),(1,\eta_2,\eta_3)$ is trivial and hence $X$ has also a base of cardinality $3$. Therefore, $X$ is not IBIS.
}
\end{example}

\begin{example}\label{example6}
{\rm Let $G$ be the Ree group $\mathrm{Ree}(q)={}^2G_2(q)$ in its $2$-transitive action of degree $q^3+1$ on the points of the Ree unital $\Omega$. There the argument is similar to Example~\ref{example5} because the stabilizer of two distinct points acts regularly on each orbit on the remaining points.

Observe that, when $q=3$, we have $\mathrm{Ree}(q)=\mathrm{P}\Gamma\mathrm{L}_2(8)$, which is not simple, and hence when
$q=3$, we also have the IBIS example $G=\mathrm{PSL}_2(8)$ of degree $3^3+1=28$.}
\end{example}

\begin{example}\label{example7}
{\rm Let $G$ be the group $\mathrm{PSL}_2(q)=\mathrm{SL}_2(q)$ with $q$ is even and let $H$ be a maximal subgroup of $G$ isomorphic to the dihedral group of order $2(q+1)$. Observe that $H$ is maximal in $G$. Let $\Omega$ be the set of conjugates of $H$ in $G$. Clearly, $G$ acts primitively on $\Omega$.

We claim that, for every $g\in G\setminus H$, $|H\cap H^g|=2$. Suppose that $H\cap H^g$ is divisible by some odd prime number $r$ and let $R$ be a Sylow $r$-subgroup of $H\cap H^g$. As $H$ and $H^g$ are dihedral groups of order $2(q+1)$, we deduce
$R\unlhd H$ and $R\unlhd H^g$. Thus $R\unlhd \langle H,H^g\rangle$. As $g\notin G\setminus H$ and as $H=\nor G H$, we get $R\unlhd G$, contradicting the fact that $G$ is a simple group. Therefore, $|H\cap H^g|\le 2$, for every $g\in G$. 
Now, let $\Sigma$ be the bipartite graph having vertex set $\Omega\cup I$, where $I$ is the set of involutions of $G$. We declare $K\in\Omega$ adjacent to $x\in I$ when $x\in K$. Now, $G$ acts as a group of automorphisms of this bipartite graph. Moreover, $G$ acts transitively on the two sides of the bipartition $\Omega\cup I$, because the involutions in $G$ form a $G$-conjugacy class. Since each element of $K\in \Omega$ contains $q+1$ involutions, we deduce that $K$ is adjacent to $q+1$ involutions.  We infer that each involution is adjacent to $$\frac{|\Omega|(q+1)}{|I|}=\frac{\frac{q}{2}(q-1)(q+1)}{q^2-1}=\frac{q}{2}$$
elements of $\Omega$. Now, let $\iota_1,\ldots,\iota_{q+1}$ be the $q+1$ involutions of $H$ and let 
$\Omega_i:=\{K\in\Omega\mid \iota_i\in K\}$. Now, for every two distinct indices $i$ and $j$, we have
$\Omega_i\cap\Omega_j=\{H\}$, because otherwise we contradict the fact that $|H\cap H^g|\le 2$, for every $g\in G\setminus H$. Therefore
$$\left|\bigcup_{i=1}^{q+1}\Omega_i\right|=(q+1)\left(\frac{q}{2}-1\right)+1=\frac{q}{2}(q-1)=|\Omega|.$$
This shows that $|H\cap H^g|=2$, for every $g\in G\setminus H$.

From the previous paragraph, it follows that, in the action of $G$ on $\Omega$, the stabilizer of any two distinct points has order $2$. It follows that $G$ is IBIS.

Arguing as in Example~\ref{example1}, it can be verified that, $G<X\le\mathrm{Aut}(G)$ is IBIS if and only if $q=2^f$ and $f$ is a prime number.}
\end{example}

\begin{example}\label{example8}{\rm Let $G:=\mathrm{P}\Gamma\mathrm{L}_2(q)$, where $q:=2^p$ for some prime number $p$ acting on $$\Omega:=\{\{V_1,V_2\}\mid V=V_1\oplus V_2, \dim V_1=\dim V_2=1\},$$
where $V=\mathbb{F}_q^2$. This action is primitive, having stabilizer a maximal subgroup in the Aschbacher class $\mathcal{C}_2$. We claim that $G$ is IBIS.

We consider the canonical basis $e_1,e_2$ of $\mathbb{F}_q^2$ giving the element $\alpha:=\{\langle e_1\rangle,\langle e_2\rangle\}$ of $\Omega$. The stabilizer of $\alpha$ is
$$
G_\alpha=\left\langle
\begin{pmatrix}a&0\\0&a^{-1}\end{pmatrix},
\begin{pmatrix}0&1\\1&0\end{pmatrix},\tau\mid a,b\in \mathbb{F}_q\setminus\{0\},\tau\in\mathrm{Gal}(\mathbb{F}_{q}/\mathbb{F}_2)
\right\rangle.$$

We now consider another direct decomposition $\beta$ of $\mathbb{F}_{q}^2$. We distinguish our analysis depending whether $\alpha\cap\beta\ne\emptyset$ or $\alpha\cap\beta=\emptyset$. 
Suppose first $\alpha\cap\beta\ne\emptyset$. Without loss of generality, we may suppose that $\langle e_1\rangle\in\beta$. Thus $\beta=\{\langle e_1\rangle,\langle xe_1+e_2\rangle\}$, for some $x\in\mathbb{F}_q\setminus\{0\}$. Now, if we take $a\in\mathbb{F}_q$ with $a^2=x^{-1}$ we obtain that
$$(xe_1+e_2)\begin{pmatrix}a&0\\0&a^{-1}\end{pmatrix}=axe_1+a^{-1}e_2=a^{-1}(a^2xe_1+e_2)=a^{-1}(e_1+e_2)\in \langle e_1+e_2\rangle .$$
Therefore we may suppose that $x=1$. Now $G_{\alpha}\cap G_{\beta}$ is the stabilizer in $G$ of the frame $\langle e_1\rangle$, $\langle e_2\rangle$, $\langle e_1+e_2\rangle$. Therefore $G_\alpha\cap G_\beta= \mathrm{Gal}(\mathbb{F}_{2^p}/\mathbb{F}_2)$ has prime order $p$.

Suppose next $\alpha\cap\beta=\emptyset$.  Thus $\beta=\{\langle e_1+xe_2\rangle,\langle ye_1+e_2\rangle\}$, for some $x,y\in\mathbb{F}_q\setminus\{0\}$ with $xy\ne 1$. Using the same trick as above, we see that we may take $x=1$ and hence $\beta=\{\langle e_1+e_2\rangle,\langle ye_1+e_2\rangle\}$, for some $y\in\mathbb{F}_q\setminus\{0,1\}$. Let $a\in\mathbb{F}_q$ with $y=a^{-2}$. An easy computation shows that
$$G_{\alpha}\cap G_{\beta}=\left\langle\begin{pmatrix}0&a\\a^{-1}&0\end{pmatrix}\right\rangle$$
has order $2$.

We have shown that, given any two distinct points $\alpha,\beta$ of $\Omega$ we have $|G_\alpha\cap G_\beta|\in \{2,p\}$ and hence $G$ is IBIS.
}
\end{example}

\begin{example}\label{example9}{\rm Let $f\ge 4$ be an even integer, let $q:=2^f$ and let $q_0:=q^{1/f}$. The group $G:=\mathrm{SL}_2(q)$ acting on the cosets of $H:=\mathrm{SL}_2(q_0)$ is IBIS. Now, we have $H={\bf C}_G(\alpha)$, where $\alpha:\mathbb{F}_q\to\mathbb{F}_q$ is the field automorphism of order $2$ defined by $\zeta^\alpha=\zeta^{q_0}$ $\forall\zeta\in\mathbb{F}_q$.

For every $g\in G\setminus H$, we have $$H\cap H^g={\bf C}_G(\alpha)\cap {\bf C}_G(\alpha^g)={\bf C}_H(\alpha^g)={\bf C}_H(\alpha g^{-1}\alpha g).$$
Set $g_\alpha:=\alpha g^{-1}\alpha g$. This shows that $H\cap H^g$ is either a Sylow $2$-subgroup of $H$ (when $g_\alpha$ has order $2$), or cyclic of order $q_0-1$ (when the order of $g_\alpha$ divides $q-1$), or cyclic of order $q_0+1$ (when the order of $g_\alpha$ divides $q+1$).

Now, let $h$ be in a $G$-coset not fixed by $H\cap H^g$ and let $h_\alpha:=\alpha h^{-1}\alpha h$. From the paragraph above, we have $H\cap H^g\cap H^h={\bf C}_H(\langle g_\alpha,h_\alpha\rangle )=1$ and hence $G$ is IBIS. Indeed, if $H\cap H^g\cap H^h\ne 1$, then $H\cap H^g={\bf C}_H(g_\alpha)={\bf C}_H(h_\alpha)=H\cap H^h$, contradicting the fact that $H\cap H^g$ does not fix the coset $Hh$.
}
\end{example}
\begin{example}\label{example10}{\rm
The group $\mathrm{Alt}(7)$ of degree $15$: this arises from the embedding $\mathrm{Alt}(7)\le \mathrm{Alt}(8)\cong\mathrm{SL}_4(2)$.
}
\end{example}

\end{document}